\documentclass[11pt]{article}
 \usepackage{graphicx,amsmath,amsfonts,amssymb,color,mathrsfs}

\usepackage{theorem}

 \newenvironment{proof}{\noindent {\bf \sc Proof.} \rm}{\QED}

\catcode`@=12

\newtheorem{thm}{\noindent Theorem}[section]
\newtheorem{lem}{\noindent Lemma}[section]

\theoremheaderfont{\normalfont\bfseries}
\theorembodyfont{\slshape} {\theorembodyfont{\rmfamily}} {\theorembodyfont{\rmfamily}
}
{\theorembodyfont{\rmfamily} } {\theorembodyfont{\rmfamily}
}
{\theorembodyfont{\rmfamily} } {\theorembodyfont{\rmfamily}
\newtheorem{rem}{\noindent Remark}[section]}
{\theorembodyfont{\rmfamily} } {\theorembodyfont{\rmfamily} } {\theorembodyfont{\rmfamily}
}

 \def\beqlb{\begin{eqnarray}}\def\eeqlb{\end{eqnarray}}
 \def\beqnn{\begin{eqnarray*}}\def\eeqnn{\end{eqnarray*}}
 \numberwithin{equation}{section}
\def\qed{\hfill$\square$\smallskip}

\makeatletter
 \newif\if@borderstar
 \def\bordermatrix{\@ifnextchar*{%
 \@borderstartrue\@bordermatrix@i}{\@borderstarfalse\@bordermatrix@i*}%
 }
 \def\@bordermatrix@i*{\@ifnextchar[{\@bordermatrix@ii}{\@bordermatrix@ii[()]}}
 \def\@bordermatrix@ii[#1]#2{%
 \begingroup
 \m@th\@tempdima8.75\p@\setbox\z@\vbox{%
 \def\cr{\crcr\noalign{\kern 2\p@\global\let\cr\endline }}%
 \ialign {$##$\hfil\kern 2\p@\kern\@tempdima & \thinspace %
 \hfil $##$\hfil && \quad\hfil $##$\hfil\crcr\omit\strut %
 \hfil\crcr\noalign{\kern -\baselineskip}#2\crcr\omit %
 \strut\cr}}%
 \setbox\tw@\vbox{\unvcopy\z@\global\setbox\@ne\lastbox}%
 \setbox\tw@\hbox{\unhbox\@ne\unskip\global\setbox\@ne\lastbox}%
 \setbox\tw@\hbox{%
 $\kern\wd\@ne\kern -\@tempdima\left\@firstoftwo#1%
 \if@borderstar\kern2pt\else\kern -\wd\@ne\fi%
 \global\setbox\@ne\vbox{\box\@ne\if@borderstar\else\kern 2\p@\fi}%
 \vcenter{\if@borderstar\else\kern -\ht\@ne\fi%
 \unvbox\z@\kern-\if@borderstar2\fi\baselineskip}%
 \if@borderstar\kern-2\@tempdima\kern2\p@\else\,\fi\right\@secondoftwo#1 $%
 }\null \;\vbox{\kern\ht\@ne\box\tw@}%
 \endgroup
 }
 \makeatother

\begin{document}
\title{Exact Tail Asymptotics --- Revisit of a Retrial Queue with Two Input Streams and Two Orbits \\
\Large (In memory of Dr. Jesus R. Artalejo)
\footnotetext{$\,$ \hspace{-6.2ex} $^{*}$ Postal address:
School of Mathematics and Statistics, Central South University,
Changsha, China, 410075
 \newline
$^{**}$ Postal address: School of Mathematics and Statistics,
Carleton University, Ottawa, ON Canada K1S 5B6 }}
\author{Yang Song$^{*}$, Zaiming Liu$^{*}$ and Yiqiang Q. Zhao$^{**}$}
\date{Revised, April 26, 2015}

\maketitle
\begin{abstract}
We revisit a single-server retrial queue with two independent Poisson streams (corresponding to two types of customers) and two orbits. The size of each orbit is infinite. The exponential server (with a rate independent of the type of customers) can hold at most one customer at a time and there is no waiting room. Upon arrival, if a type $i$ customer $(i=1,2)$ finds a busy server, it will join the type $i$ orbit. After an exponential time with a constant (retrial) rate $\mu_i$, an type $i$ customer attempts to get service. This model has been recently studied by Avrachenkov, Nain and Yechiali~\cite{ANY2014} by solving a Riemann-Hilbert boundary value problem. One may notice that, this model is not a random walk in the quarter plane. Instead, it can be viewed as a random walk in the quarter plane modulated by a two-state Markov chain, or a two-dimensional quasi-birth-and-death (QBD) process. The special structure of this chain allows us to deal with the fundamental form corresponding to one state of the chain at a time, and therefore it can be studied through a boundary value problem. Inspired by this fact, in this paper, we focus on the tail asymptotic behaviour of the stationary joint probability distribution of the two orbits with either an idle or busy server by using the kernel method, a different one that does not require a full determination of the unknown generating function. To take advantage of existing literature results on the kernel method, we identify a censored random walk, which is an usual walk in the quarter plane.
This technique can also be used for other random walks modulated by a finite-state Markov chain with a similar structure property.
\medskip

\medskip


\end{abstract}

\section{Introduction}

In this paper, we revisit a single server retrial queue system with two orbits and no waiting room, which has been studied by Avrachenkov, Nain and Yechiali~\cite{ANY2014}. The analysis in \cite{ANY2014} is based on the solution of a Riemann-Hilbert boundary value problem, while our focus is on exact tail asymptotics for the joint stationary distribution of the two orbits under a busy or idle state of the server, using the kernel method, a different method that does not require a full determination of the unknown generating function.  In this system, there are two independent exogenous Poisson streams (representing two types of customers) flowing into the server, and the server can hold at most one customer at a time. Upon arrival, if a type $i$ customer finds a busy server, it will join its orbit and wait for retrial at a specified exponential rate for the customers of type $i$. Such a queueing system could serve as a model for two competing job streams in a carrier sensing multiple access system, and it has an application in a local area computer network (LAN) as explained in \cite{ANY2014}.

Retrial queueing systems have been attracting researchers' attention for many years (e.g., \cite{AG2008,A2010,F1990,YT1987} and references therein). Much of the previous work lays the emphasis on performance measures, such as the mean size of the orbit, the average number of the customers in the system, the average waiting time among others. We also notice that stationary tail asymptotic analysis has recently become one of the central research topics for retrial queues due to not only its own importance, but also its applications in approximation and performance bounds. For example, in \cite{SLL2006}, Shang, Liu and Li proved that the stationary queue length of the $M/G/1$ retrial queue has a subexponential tail if the queue length of the corresponding $M/G/1$ queue has a tail of the same type. Kim, Kim and Kim extended the study on the $M/G/1$ retrial queue in \cite{KKK2007} by Kim, Kim and Ko to a $MAP/G/1$ retrial queue, and obtained tail asymptotics for the queue size distribution in \cite{KKK2010}. By adopting matrix-analytic theory and the censoring technique in \cite{LWZ2012}, Liu, Wang and Zhao studied the $M/M/c$ retrial queues with non-persistent customers and obtained tail asymptotics for the joint stationary distribution of the number of retrial customers in the orbit and the number of busy servers.

Most of the studies on retrial queues assumed a single type of customers flowing into the system, and references on retrial systems with multi-class customers are quite limited. The model studied by Avrachenkov, Nain and Yechiali in \cite{ANY2014} and again in this paper is such a system. This model is an example of the two-dimensional QBD process (for example, see Ozawa~\cite{Ozawa:2013}), or the random walk in the quarter plane modulated by a two-state Markov chain (another example of retrial queues having this structure is Li and Zhao~\cite{Li-Zhao:2005}). In \cite{ANY2014}, the authors showed how this modulated model is converted to a scalar fundamental form, which can be solved in terms of a Riemann-Hilbert boundary value problem (BVP) due to its special structure of this system.
Motivated by this, we extend their research on this model by considering the tail asymptotic behaviour of the stationary joint probability distribution of the two orbits with either an idle or a busy server, by using the kernel method --- a different method that does not require a full determination of the unknown generating function.  For more details about the kernel method,  readers may refer to \cite{FIM1999,LZ2009,LZ2011,LZ2012,LTZ2013}. We point out that tail asymptotic properties for Markov modulated or more general block-structured random walks have also been studied by using other methods, for example in \cite{Miyazawa-Zhao:2004, SMZ:2006, KMZ:2010, Miyazawa:2015}.

The main contribution in this paper includes: (1) the characterization of the tail asymptotic properties in the joint distribution for a large queue $i$ ($i=1, 2$) with either an idle or a busy server. A total of three types of properties are identified (see Theorems~\ref{theorem:6.2}, \ref{theorem:6.3}, and \ref{theorem:6.4} for the case of a busy server, and Theorems~\ref{theorem:7.1} and \ref{theorem:7.2} for the case of a idle server); and (2) an illustration on how to convert a matrix-form fundamental form for the Markov modulated random walk into a (scalar)
functional form corresponding to one state of the chain, through a censored Markov chain or solving the matrix-form fundamental form (see remarks in the last section). Therefore, it can be studied by the kernel method.

The rest of the paper is organized as follows: Section~2 provides the model description; Section~3 identifies the censored random walk in the quarter plane; dominant singularities of the unknown generating function are located in Section~4, while the detailed asymptotic property of the unknown function at its dominant singularity is discussed in Section~5; exact tail asymptotic properties for stationary probabilities of the system, which are our main results, and their detail proofs are presented in Section~6. Concluding remarks are made in the final section.

\section{Model description}

In this paper, we consider a single server queueing system with two independent Poisson streams of arrivals and two retrial orbits, the same system studied in \cite{ANY2014}. Following \cite{ANY2014}, the two arrival rates are denoted by $\lambda_i$, $i=1,2$, with $\lambda=\lambda_1+\lambda_2$. The server can hold at most one customer at a time (without a waiting room). It means that when the server is busy, an arriving type $i$ customer will join in orbit $i$ of infinity capacity. Retrials from all customers in orbit $i$ for service are characterized by a Poisson process with constant rate $\mu_i$. The service time for each customer is independent of its type and follows an exponential distribution with rate $\mu$.  The retrial mechanism imposed can be a model when only the customer at head of the line (orbit) is allowed for retrial.

Let $I(t)$ be the state of the server (either idle or busy), or the number of customers in the server, and let $Q_i(t)$ denote the number of customers in orbit $i$ at time $t$ for $i=1,2$. Then, it generates a continuous time Markov chain $X(t)=\{(Q_1(t), Q_2(t), I(t)):t\in [0,\infty]\}$ on the state space $\{0,1,\ldots\}\times\{0,1,\ldots\}\times\{0,1\}$. From Avrachenkov, Nain and Yechiali \cite{ANY2014}, we know that the system is stable if and only if $\lambda(\lambda_1+\mu_1)<\mu\mu_1$ and $\lambda(\lambda_2+\mu_2)<\mu\mu_2$.
Under the stability condition, the unique stationary probability vector for the system is denoted by
$\Pi_{m,n}=(\pi_{m,n}(0),\pi_{m,n}(1))$
for $m,n=0,1, \ldots$. For the purpose of finding the stationary distribution, we consider the corresponding discrete time Markov chain through the uniformization technique. Without loss of generality, we assume that $\lambda+\mu+\mu_1+\mu_2=1$. For this discrete time chain, a transition diagram, partitioned according to the state of the server, is depicted in Figure~1, where
\begin{equation*}
{A_{1,0}=A_{1,0}^{(1)}=A_{1,0}^{(2)}=A_{1,0}^{(0)}}=
\left(
\begin{array}{cccc}
0 & 0 \\
0 & \lambda_1
\end{array}\right),
\end{equation*}

\begin{equation*}
{A_{0,1}=A_{0,1}^{(1)}=A_{0,1}^{(2)}=A_{0,1}^{(0)}}=
\left(
\begin{array}{cccc}
0 & 0 \\
0 & \lambda_2
\end{array}\right),
\end{equation*}

\begin{equation*}
{A_{-1,0}=A_{-1,0}^{(1)}}=
\left(
\begin{array}{cccc}
0 & \mu_1 \\
0 & 0
\end{array}\right),
\end{equation*}

\begin{equation*}
{A_{0,-1}=A_{0,-1}^{(2)}}=
\left(
\begin{array}{cccc}
0 & \mu_2 \\
0 & 0
\end{array}\right),
\end{equation*}

\begin{equation*}
{A_{0,0}}=
\left(
\begin{array}{cccc}
\mu & \lambda \\
\mu & \mu_1+\mu_2
\end{array}\right),
\end{equation*}

\begin{equation*}
{A_{0,0}^{(1)}}=
\left(
\begin{array}{cccc}
\mu+\mu_2 & \lambda \\
\mu & \mu_1+\mu_2
\end{array}\right),
\end{equation*}

\begin{equation*}
{A_{0,0}^{(2)}}=
\left(
\begin{array}{cccc}
\mu+\mu_1 & \lambda \\
\mu & \mu_1+\mu_2
\end{array}\right),
\end{equation*}

\begin{equation*}
{A_{0,0}^{(0)}}=
\left(
\begin{array}{cccc}
\mu+\mu_1+\mu_2 & \lambda \\
\mu & \mu_1+\mu_2
\end{array}\right).
\end{equation*}
\begin{center}\includegraphics[width=8cm,height=7cm]{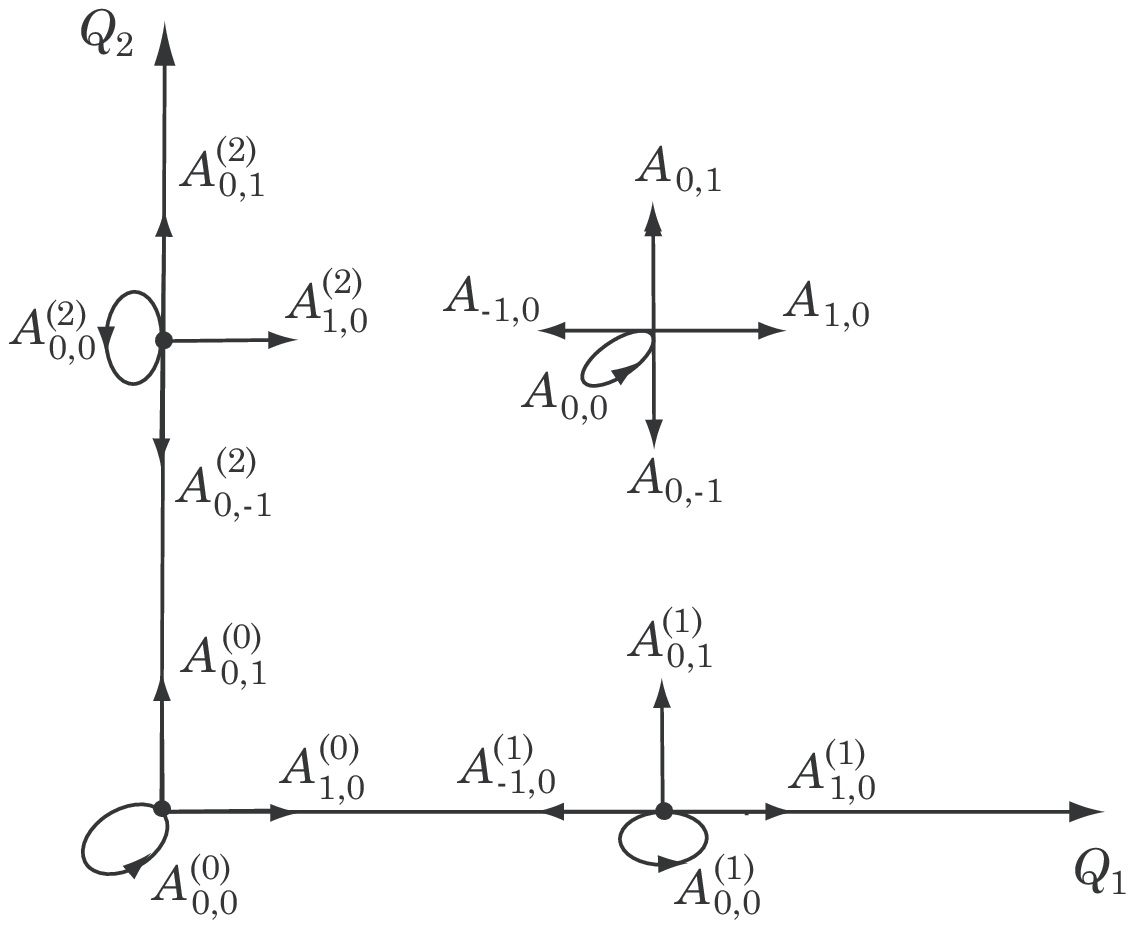}\\
{\footnotesize {\bf Figure 1} Matrix transition diagram }
\end{center}

We define the probability generating function (PGF) $P^{(k)}(x,y)$ for the stationary probabilities $\pi_{m,n}(k)$ as
\begin{equation*}
P^{(k)}(x,y)=\sum_{m=0}^{\infty}\sum_{n=0}^{\infty}\pi_{m,n}(k)x^{m}y^{n},\quad|x|\leq1,|y|\leq1,\quad k=0,1,
\end{equation*}
and denote
\begin{equation*}
P(x,y)=(P^{(0)}(x,y),P^{(1)}(x,y)).
\end{equation*}
Following the idea in Fayolle, Iasnogorodski and Malyshev \cite{FIM1999}, we can obtain the (matrix-form) fundamental form for the Markov modulated random walk in the quarter plane:
\beqlb\label{2-1}
P(x,y)H(x,y)=P(x,0)H_{1}(x,y)+P(0,y)H_{2}(x,y)+\Pi_{0,0}H_{0}(x,y),
\eeqlb
where
\beqnn
H(x,y)&=&-\bar{h}(x,y),
\\H_1(x,y)&=&-\bar{h}(x,y)+\bar{h}_1(x,y)y,
\\H_2(x,y)&=&-\bar{h}(x,y)+\bar{h}_2(x,y)x,
\\H_0(x,y)&=&\bar{h}_0(x,y)xy+\bar{h}(x,y)-\bar{h}_1(x,y)y-\bar{h}_2(x,y)x
\eeqnn
with
\beqnn
\bar{h}(x,y)&=&xy\Big(\sum_{i=-1}^1\sum_{j=-1}^1 A_{i,j}x^{i}y^{j}-I\Big),\quad
\\ \bar{h}_1(x,y&)=&x\Big(\sum_{i=-1}^1\sum_{j=0}^1A^{(1)}_{i,j}x^{i}y^{j}-I\Big),\quad
\\ \bar{h}_2(x,y)&=&y\Big(\sum_{i=0}^1\sum_{j=-1}^1A^{(2)}_{i,j}x^{i}y^{j}-I\Big),\quad
\\ \bar{h}_0(x,y)&=&\sum_{i=0}^1\sum_{j=0}^1A^{(0)}_{i,j}x^{i}y^{j}-I.
\eeqnn
For detailed derivation, see the work in \cite{G-thesis}.

\begin{rem}
It is worthwhile to mention that the fundamental form (1.3.6) in \cite{FIM1999} is for the generating function excluding boundary probabilities, while ours is for the complete joint probability vector. In \cite{G-thesis}, it is pointed out that these two forms (for Markov modulated random walks) are equivalent. In fact, for $k=0,1$ let
\[
    \pi^{(k)}(x,y)=\sum_{m=1}^{\infty}\sum_{n=1}^{\infty}\pi_{m,n}(k)x^{m-1}y^{n-1},
\]
and $\pi(x,y)=(\pi^{(0)}(x,y), \pi^{(1)}(x,y))$, then
\[
    -\pi(x,y) \bar{h}(x,y) = \pi(x,0) \bar{h}_1(x,y) + \pi(0,y) \bar{h}_2(x,y) + \Pi_{0,0} \bar{h}_0(x,y)
\]
by noticing that
\begin{eqnarray*}
   \pi_1^{(k)}(x) = \pi^{(k)}(x,0)&=&\sum_{m=1}^{\infty}\pi_{m,0}(k)x^{m-1},
\\ \pi_2^{(k)}(y) = \pi^{(k)}(0,y)&=&\sum_{n=1}^{\infty}\pi_{0,n}(k)y^{n-1}.
\end{eqnarray*}
\end{rem}

For the retrial queueing system studied in this paper, after some calculations, we have
\begin{equation*}
H(x,y)=
\left(
\begin{array}{cccc}
(\lambda+\mu_1+\mu_2)xy & -(\mu_2x+\mu_1y+\lambda xy)\\
-\mu xy & -[\lambda_2 xy^2+\lambda_1 x^2 y-(\lambda+\mu)xy]
\end{array}\right),
\end{equation*}

\begin{equation*}
H_1(x,y)=
\left(
\begin{array}{cccc}
\mu_2 xy & -\mu_2 x\\
0 & 0
\end{array}\right),
\end{equation*}
\begin{equation*}
H_2(x,y)=
\left(
\begin{array}{cccc}
\mu_1 xy & -\mu_1 y\\
0 & 0
\end{array}\right),
\end{equation*}
\begin{equation*}
H_0(x,y)=
\left(
\begin{array}{cccc}
0 & 0\\
0 & 0
\end{array}\right).
\end{equation*}
Hence, the fundamental form (\ref{2-1}) can be simplified as
\begin{equation*} 
P(x,y)H(x,y)=P(x,0)H_{1}(x,y)+P(0,y)H_{2}(x,y).
\end{equation*}
Equivalently,
\begin{align*}
(P^{(0)}(x,y),P^{(1)}(x,y))\left(
\begin{array}{cccc}
(\lambda+\mu_1+\mu_2)xy & -(\mu_2x+\mu_1y+\lambda xy)\\
-\mu xy & -[\lambda_2 xy^2+\lambda_1 x^2 y-(\lambda+\mu)xy]
\end{array}\right)
\\=(P^{(0)}(x,0),P^{(1)}(x,0))\left(
\begin{array}{cccc}
\mu_2 xy & -\mu_2 x\\
0 & 0
\end{array}\right)+(P^{(0)}(0,y),P^{(1)}(0,y))\left(
\begin{array}{cccc}
\mu_1 xy & -\mu_1 y\\
0 & 0
\end{array}\right),
\end{align*}
or,
\begin{flushleft}
\begin{equation}\label{2-4}
(\lambda+\mu_1+\mu_2)P^{(0)}(x,y)-\mu P^{(1)}(x,y)=\mu_2 P^{(0)}(x,0)+\mu_1P^{(0)}(0,y),
\end{equation}
\begin{equation}\label{2-5}
(\lambda xy+\mu_1 y+\mu_2 x)P^{(0)}(x,y)+[\lambda_1 x+\lambda_2 y-(\lambda+\mu)]xy P^{(1)}(x,y)=\mu_2 xP^{(0)}(x,0)+\mu_1 y P^{(0)}(0,y).
\end{equation}
\end{flushleft}
\eqref{2-4} and \eqref{2-5} are identical to equations (18) and (19) in \cite{ANY2014}, derived from direct calculations.

\section{Censored Markov chain}

One may notice that equations (\ref{2-4}) and (\ref{2-5}) provide a relationship between generating functions for an idle server and for a busy server.
Therefore, we start our analysis for a busy server since in this case, the censored Markov chain can be expressed explicitly. This censored Markov chain is a random walk in the quarter plane, which has been extensively studied in the literature. To this end, we first consider the uniformized discrete time Markov chain of the continuous time chain $X(t)$ for the retrial model with uniformization parameter $\theta = \lambda + \mu +\mu_1 +\mu_2=1$.  We partition the transition matrix $P$ of the uniformized chain according to the server state and then consider the censored chain to the set of states of a busy server. Specifically, let
 $X(n)=\{(Q_n^{(1)},Q_n^{(2)}, I_n)\}$ be the uniformized chain on the state space $\{0, 1, \ldots \} \times \{0, 1, \ldots \} \times \{0, 1\}$ and let
$E=\{0,1,\ldots\}\times\{0,1,\ldots\}\times\{1\}$ and $E^c=\{0,1,\ldots\}\times\{0,1,\ldots\}\times\{0\}$. Partition the transition matrix $P$ according to $E$ and its complement $E^c$ into:
\[
    P = \bordermatrix[{()}]{
    & \text{$E^c$} & \text{$E$} \cr
\text{$E^c$} & P_{00}& P_{01}   \cr
\text{$E$} &  P_{10} & P_{11}
},
\]
where using the lexicographical order for states of $(Q_n^{(1)},Q_n^{(2)})$, $P_{ij}$ can be expressed as
\[
P_{00}=\left(
    \begin{array}{cccccc}
     &A_{0} &  \\
     &      &A_{1} & \\
     &      &      &A_{1} & \\
     &      &      &      &\ddots
    \end{array}
  \right),\quad\quad
P_{01}=\left(
    \begin{array}{cccccc}
     &B_{0} &  \\
     &B_{1} &B_{0} & \\
     &      &B_{1} &B_{0}  & \\
     &      &      &\ddots &\ddots
    \end{array}
  \right),
\]

\[
P_{10}=\left(
    \begin{array}{cccccc}
     &C_0 &  \\
     &          &C_0  &\\
     &          &       &\ddots
    \end{array}
  \right);\quad\quad
 P_{11}=\left(
    \begin{array}{cccccc}
     &D_{0} &D_{1}  \\
     &      &D_{0} &D_{1} \\
     &      &      &\ddots &\ddots
    \end{array}
  \right),
\]
with
\[
A_{0}=\left(
    \begin{array}{cccccc}
     &\mu+\mu_1+\mu_2 &  \\
     &                &\mu+\mu_1 & \\
     &                &          &\mu+\mu_1 & \\
     &                &          &          &\ddots
    \end{array}
  \right),\quad\quad
A_{1}=\left(
    \begin{array}{cccccc}
     &\mu+\mu_2 &  \\
     &          &\mu  &\\
     &          &     &\mu  &  \\
     &          &     &     & \ddots
    \end{array}
  \right),
\]

\[
B_{0}=\left(
    \begin{array}{cccccc}
     &\lambda &  \\
     &\mu_2   &\lambda & \\
     &        &\mu_2   &\lambda  & \\
     &        &        &\ddots &\ddots
    \end{array}
  \right),\quad\quad
B_{1}=\left(
    \begin{array}{cccccc}
     &\mu_1 &  \\
     &          &\mu_1  &\\
     &          &       &\ddots
    \end{array}
  \right), \quad\quad
C_{0}=\left(
    \begin{array}{cccccc}
     &\mu       &  \\
     &          &\mu    &  \\
     &          &       &\ddots
    \end{array}
  \right),
\]
\[
D_{0}=\left(
    \begin{array}{cccccc}
     &\mu_1+\mu_2  &\lambda_2  \\
     &             &\mu_1+\mu_2 &\lambda_2  \\
     &             &            &\ddots   &\ddots
    \end{array}
  \right),\quad\quad
D_{1}=\left(
    \begin{array}{cccccc}
     &\lambda_1    &  \\
     &             &\lambda_1    &  \\
     &             &             & \ddots
    \end{array}
  \right).
\]

Notice that $P_{00}$ is diagonal, it is straightforward to have the fundamental matrix of $P_{00}$ as follows:
\[
    \hat{P}_{00}=\sum_{n=0}^{\infty} P_{00}^n = \text{diag} (\hat{A}_0, \hat{A}_1, \hat{A}_1, \ldots),
\]
where
\begin{align*}
    \hat{A}_0 &= \text{diag} \left ( \frac{1}{\lambda}, \frac{1}{\lambda+\mu_2}, \frac{1}{\lambda+\mu_2}, \ldots \right), \\
    \hat{A}_1 &= \text{diag} \left ( \frac{1}{\lambda+\mu_1}, \frac{1}{\lambda+\mu_1+\mu_2}, \frac{1}{\lambda+\mu_1+\mu_2}, \ldots \right).
\end{align*}
Furthermore, notice that $P_{10}$ is also diagonal, and therefore the censored chain to $E$ can be easily computed as
\[
    P^{(E)} = P_{11}+P_{10}\hat{P}_{00}P_{01} =  \left(
    \begin{array}{cccccc}
     D_0 + \mu\hat{A}_0B_0 & D_1    &  \\
      \mu\hat{A}_1B_1 & D_0 + \mu\hat{A}_1B_0  & D_1  &  \\
     & \mu\hat{A}_1B_1 & D_0 + \mu\hat{A}_1B_0  & D_1  &  \\
     &  & \ddots   &  \ddots   & \ddots &
    \end{array}
  \right).
\]
This censored chain is an example, referred to the simple random walk, of the random walks in the quarter plane studied in \cite{FIM1999}, whose transition diagram is depicted in Figure~2.

\begin{center}\includegraphics[width=8cm,height=7cm]{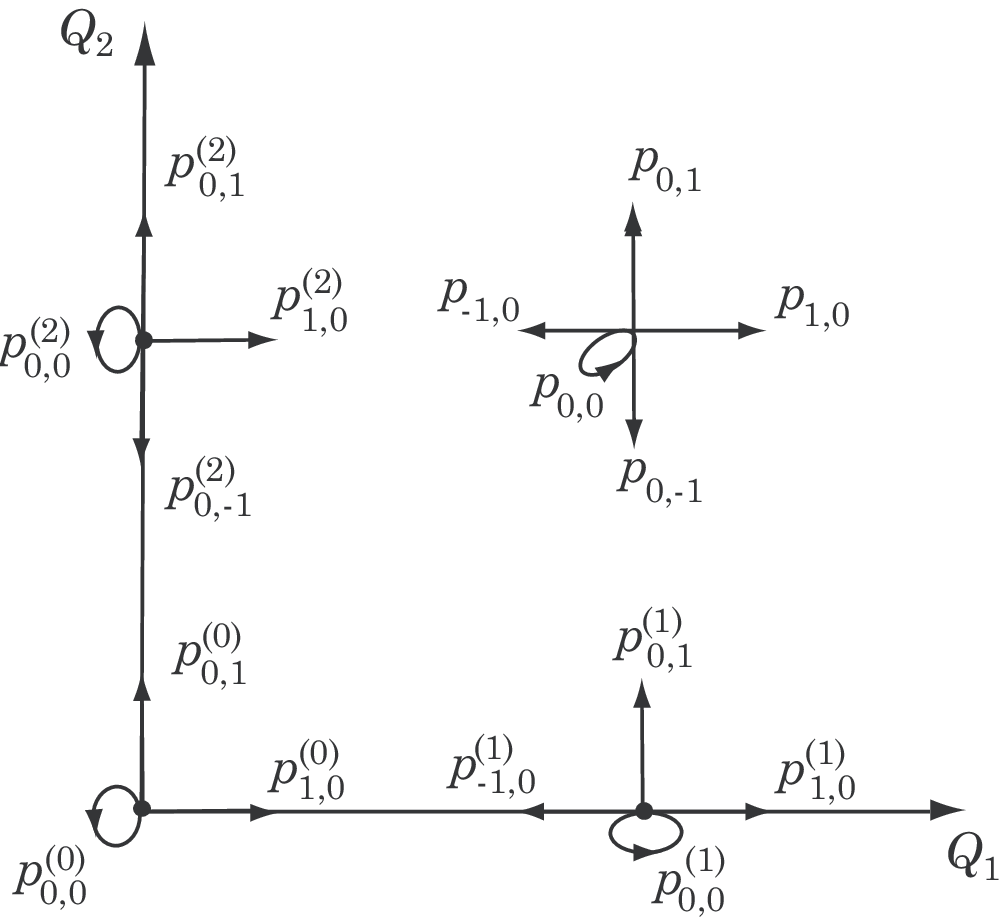}\\
{\footnotesize {\bf Figure 2} Transition diagram of the censored random walk}
\end{center}

In our case,
\beqnn
p_{1,0}=p_{1,0}^{(1)}=p_{1,0}^{(2)}=p_{1,0}^{(0)}=\lambda_1, ~~ p_{0,1}=p_{0,1}^{(1)}=p_{0,1}^{(2)}=p_{0,1}^{(0)}=\lambda_2;
\eeqnn
\beqnn
p_{-1,0}=\frac{\hat{\mu}_1}{\lambda+\mu_1+\mu_2}, ~~ p_{0,-1}=\frac{\hat{\mu}_2}{\lambda+\mu_1+\mu_2}, ~~ p_{0,0}=\mu_1+\mu_2+\frac{\lambda\mu}{\lambda+\mu_1+\mu_2};
\eeqnn
\beqnn
p_{-1,0}^{(1)}=\frac{\hat{\mu}_1}{\lambda+\mu_1}, ~~ p_{0,0}^{(1)}=\mu_1+\mu_2+\frac{\lambda\mu}{\lambda+\mu_1};
\eeqnn
\beqnn
p_{0,-1}^{(2)}=\frac{\hat{\mu}_2}{\lambda+\mu_2}, ~~ p_{0,0}^{(2)}=\mu_1+\mu_2+\frac{\lambda\mu}{\lambda+\mu_2};
\eeqnn
\beqnn
p_{0,0}^{(0)}=\mu_1+\mu_2+\mu,
\eeqnn
where $\hat{\mu}_i=\mu\mu_i$ for $i=1,2$.

Let $\alpha=\lambda+\mu_1+\mu_2=1-\mu$ and $\hat{\lambda}_i=\alpha\lambda_i$ for $i=1,2$.
It follows from \cite{ANY2014} that under the system stability condition (for the retrial queue), at least one of $\hat{\lambda}_1<\hat{\mu}_1$ and $\hat{\lambda}_2<\hat{\mu}_2$ holds. Without loss of generality, we assume that $\hat{\lambda}_1<\hat{\mu}_1$ throughout the paper. For this censored random walk, the fundamental form (equation (1.3.6) in \cite{FIM1999}) is given by:
\beqlb\label{3-15}
    -h(x,y) \pi^{(1)}(x,y)= h_1(x,y) \pi_1^{(1)}(x) + h_2(x,y) \pi_2^{(1)}(y) + h_0(x,y) \pi_{0,0}(1),
\eeqlb
where
\begin{eqnarray}\label{3-2}
     h(x,y)&=&xy\Big(\sum_{i=-1}^1\sum_{j=-1}^1 p_{i,j}x^{i}y^{j}-1\Big) = a(x)y^2+b(x)y+c(x), \\
\label{3-3}
    h_1(x,y)&=&x\Big(\sum_{i=-1}^1\sum_{j=0}^1 p^{(1)}_{i,j}x^{i}y^{j}-1\Big) = a_1(x)y+b_1(x), \\
\label{3-4}
     h_2(x,y)&=&y\Big(\sum_{i=0}^1\sum_{j=-1}^1 p^{(2)}_{i,j}x^{i}y^{j}-1\Big) =a_2(x)y^2+b_2(x)y+c_2(x), \\
\label{3-5}
    h_0(x,y)&=&\sum_{i=0}^1\sum_{j=0}^1 p^{(0)}_{i,j}x^{i}y^{j}-1 =a_0(x)y+b_0(x),
\end{eqnarray}
with
\begin{eqnarray*}
a(x)&=&p_{0,1}x, \quad b(x)=p_{-1,0}-(1-p_{0,0})x+p_{1,0}x^2, \quad c(x)=p_{0,-1}x;
\\a_1(x)&=&p_{0,1}^{(1)}x, \quad b_1(x)=p_{-1,0}^{(1)}-(1-p_{0,0}^{(1)})x+p_{1,0}^{(1)}x^2;
\\a_2(x)&=&p_{0,1}^{(2)}, \quad b_2(x)=p_{0,0}^{(2)}-1+ p_{1,0}^{(2)}x, \quad c_2(x)=p_{0,-1}^{(2)};
\\a_0(x)&=&p_{0,1}^{(0)}, \quad b_0(x)=p_{1,0}^{(0)}x-(1-p_{0,0}^{(0)}).
\end{eqnarray*}

In \cite{LZ2011, LZ2012}, a kernel method has been promoted for studying exact tail asymptotic properties for random walks in the quarter plane. In the following, we apply this method to the retrial queue model to explicitly (in terms of system parameters) characterize regions on which different tail asymptotic properties hold. First, based on the fundamental form in (\ref{3-15}), asymptotic properties at the dominant singularity for the generating function $\pi_1^{(1)}(x)$ or $\pi_2^{(1)}(y)$ are obtained, based on which regions of different exact tail asymptotic properties for probabilities $\pi_{m,n}(1)$ with a fixed value of $n$ or $m$ are identified through a Tauberian-like theorem (Theorem~\ref{thm6-1}). Then, based on the relationship given in (\ref{2-4}), the generating functions $\pi_1^{(0)}(x)$ and $\pi_2^{(0)}(y)$ are analyzed, and characterization of the exact tail asymptotic properties
for $\pi_{m,n}(0)$ is provided.

\section{Dominant singularity of $\pi_1^{(1)}(x)$}

Since discussions for dominant singularities of the two functions $\pi_1^{(1)}(x)$ and $\pi_2^{(1)}(y)$ are repetitive, we only provide details for $\pi_1^{(1)}(x)$.

According to \cite{LZ2011,LZ2012}, the dominant singularity of $\pi_1^{(1)}(x)$ is either a branch point of the Riemann surface defined by the kernel equation $h(x,y)=0$, or a pole of the function $\pi_1^{(1)}(x)$. The following two subsections are devoted to these two cases, respectively.

\subsection{Branch points for kernel equation $h(x,y)=0$}

For the censored random walk, we consider the kernel equation $h(x,y)=0$ defined by the kernel function $h(x,y)$ in the fundamental form \eqref{3-15}. Write $\alpha h(x,y)$ as a quadratic form in $y$ with coefficients that are polynomials in $x$:
\begin{equation}\label{4-1}
\alpha h(x,y)=(\hat{\lambda}_2x)y^2+\left[\hat{\lambda}_1 x^2-(\hat{\lambda}_1+\hat{\lambda}_2+\hat{\mu}_1+\hat{\mu}_2)x+\hat{\mu}_1\right]y+\hat{\mu}_2 x.
\end{equation}
For a fixed $x$, $h(x,y)=0$ has two solutions
\begin{equation*}
Y_{\pm}(x)=\frac{-\hat{b}(x)\pm\sqrt{\Delta(x)}}{2\hat{\lambda}_2 x},
\end{equation*}
where $\hat{b}(x)=\alpha b(x)=\hat{\lambda}_1 x^2-(\hat{\lambda}_1+\hat{\lambda}_2+\hat{\mu}_1+\hat{\mu}_2)x+\hat{\mu}_1$ and $\Delta(x)=b_+(x)b_-(x)$ with
\beqnn \label{4-2}
b_+(x)=\hat{b}(x)+2x\sqrt{\hat{\lambda}_2\hat{\mu}_2}=(x-1)(\hat{\lambda}_1 x-\hat{\mu}_1)-(\sqrt{\hat{\lambda}_2}-\sqrt{\hat{\mu}_2})^2x,
\eeqnn
\begin{equation*}
b_-(x)=\hat{b}(x)-2x\sqrt{\hat{\lambda}_2\hat{\mu}_2}=(x-1)(\hat{\lambda}_1 x-\hat{\mu}_1)-(\sqrt{\hat{\lambda}_2}+\sqrt{\hat{\mu}_2})^2x.
\end{equation*}

Denote the branch points by $x_i,i=1,2,3,4$, which are the zeros of $\Delta(x)$, then we have
\begin{equation}\label{4-3}
b_+(x)=\hat{\lambda}_1(x-x_2)(x-x_3)\quad \text{and} \quad b_-(x)=\hat{\lambda}_1(x-x_1)(x-x_4),
\end{equation}
where
\[
    0<x_1<x_2<1<\hat{\mu}_1/\hat{\lambda}_1\leq x_3<x_4<+\infty
\]
and
\begin{equation} \label{eqn:x3}
    x_3 = \frac{\big(\hat{\lambda}_1 + \hat{\mu}_1 \big  ) + \big (\sqrt{\hat{\lambda}_2}-\sqrt{\hat{\mu}_2} \big )^2 - \sqrt{ \Big [\big (\hat{\lambda}_1 + \hat{\mu}_1 \big ) + \big (\sqrt{\hat{\lambda}_2}-\sqrt{\hat{\mu}_2} \big )^2 \Big ]^2 - 4 \hat{\lambda}_1 \hat{\mu}_1}}{2 \hat{\lambda}_1}
\end{equation}
is a candidate of the dominant singularity of $\pi_1^{(1)}(x)$.

Consider the following cut planes:
$$\widetilde{\mathbb{C}}_x=\mathbb{C}_x-[x_3,x_4],$$
$$\widetilde{\mathbb{C}}_y=\mathbb{C}_y-[y_3,y_4],$$
$$\widetilde{\widetilde{\mathbb{C}}}_x=\mathbb{C}_x-[x_3,x_4]\cup[x_1,x_2],$$
$$\widetilde{\widetilde{\mathbb{C}}}_y=\mathbb{C}_y-[y_3,y_4]\cup[y_1,y_2],$$
where
${\mathbb{C}}_x$ and ${\mathbb{C}}_y$ are the complex planes of $x$ and $y$, respectively.
In the cut plane $\widetilde{\widetilde{\mathbb{C}}}_x$, define the two branches of $Y(x)$ by
\begin{align*}
Y_0(x)=Y_-(x) \quad \textrm{and} \quad Y_1(x)=Y_+(x) \quad if\quad |Y_-(x)|\leq|Y_+(x)| ,
\\Y_0(x)=Y_+(x) \quad \textrm{and} \quad Y_1(x)=Y_-(x) \quad if\quad |Y_-(x)|>|Y_+(x)| .
\end{align*}

Symmetrically, when $x$ and $y$ are interchanged, we also have branch points $y_i,i=1,2,3,4$, satisfying
$$0<y_1<y_2<1<y_3<y_4<+\infty$$
as well as the two branches $X_0(y)$ and $X_1(y)$ defined in a similar fashion.

Detailed properties of the branches $Y_0(x)$ and $Y_1(x)$ ($X_0(y)$ and $X_1(y)$) are needed in the asymptotic analysis for functions $\pi_1^{(1)}(x)$ and $\pi_1^{(0)}(x)$ ($\pi_2^{(1)}(y)$
and $\pi_2^{(0)}(y)$), which are presented in the following two lemmas.

\begin{lem}\label{lem4-1}
The functions $Y_i(x)$, $i=0,1,$ are meromorphic in the cut plane $\widetilde{\widetilde{\mathbb{C}}}_x$. In addition,
\begin{description}
\item[(i)] $Y_0(x)$ has one zero and no poles and $Y_1(x)$ has two poles and no zeros. Hence, $Y_0(x)$ is analytic in $\widetilde{\widetilde{\mathbb{C}}}_x$.
\item[(ii)] $|Y_0(x)|<|Y_1(x)|$ in the whole cut complex plane $\widetilde{\widetilde{\mathbb{C}}}_x$. $|Y_0(x)|=|Y_1(x)|$ takes place only on the cuts.
\item[(iii)] $|Y_0(x)|<1$ if $|x|=1$, $x\neq1$, and $Y_0(1)=\min\left(1,\frac{\hat{\mu}_2}{\hat{\lambda}_2}\right)\leq1$.
\item[(iv)] For all $x\in\mathbb{C}_x$, $|Y_0(x)|\leq\sqrt{\frac{\hat{\mu}_2}{\hat{\lambda}_2}}$ and $|Y_1(x)|\geq\sqrt{\frac{\hat{\mu}_2}{\hat{\lambda}_2}}$.
\item[(v)] If $x\in[x_1,x_2]$, then $|Y_0(x)|=\sqrt{\frac{\hat{\mu}_2}{\hat{\lambda}_2}}$ and $X_0(Y_0(x))=x$.
\end{description}
Moreover,
\begin{description}
\item[(vi)] $0<Y_0(x)\leq1$ for $1\leq x\leq\frac{\hat{\mu}_1}{\hat{\lambda}_1}$ (recall that $\hat{\lambda}_1<\hat{\mu}_1$).
\end{description}
Parallel results for $X_i(y)$, $i=0,1$ can be stated as well.
\end{lem}

\proof
See \cite{FI1979}, \cite{LZ2011} and \cite{LZ2012} for proofs of (i)--(v). Here we only detail the proof to (vi).

For $1\leq x\leq\frac{\hat{\mu}_1}{\hat{\lambda}_1}$, let $\frac{\alpha h(x,y)}{x}=0$, which leads to
$$\tilde{b}(y)+\left(\hat{\lambda}_1 x+\frac{\hat{\mu}_1}{x}\right)y=0,$$
where $\tilde{b}(y)=\hat{\lambda}_2y^2-(\hat{\lambda}_1+\hat{\lambda}_2+\hat{\mu}_1+\hat{\mu}_2)y+\hat{\mu}_2$.
Since $\hat{\lambda}_1 x+\frac{\hat{\mu}_1}{x}$ is decreasing on $\left[1,\sqrt{\frac{\hat{\mu}_1}{\hat{\lambda}_1}}\right]$ and increasing on $\left(\sqrt{\frac{\hat{\mu}_1}{\hat{\lambda}_1}},\frac{\hat{\mu}_1}{\hat{\lambda}_1}\right]$, therefore,
$2\sqrt{\hat{\lambda}_1\hat{\mu}_1}\leq\hat{\lambda}_1 x+\frac{\hat{\mu}_1}{x}\leq\hat{\lambda}_1+\hat{\mu}_1$.
For $y<0$, the inequalities
\begin{equation*}
\left\{\begin{array}{ll}
\tilde{b}(y)+2 y\sqrt{\hat{\lambda}_1\hat{\mu}_1}\geq0
\\\tilde{b}(y)+(\hat{\lambda}_1+\hat{\mu}_1)y\leq0
\end{array}
\right.
\end{equation*}
have no solutions.
For $y\geq0$, solve the following inequalities
\begin{equation*}
\left\{\begin{array}{ll}
\tilde{b}(y)+2 y\sqrt{\hat{\lambda}_1\hat{\mu}_1}\leq0
\\\tilde{b}(y)+(\hat{\lambda}_1+\hat{\mu}_1)y\geq0
\end{array}
\right.
\end{equation*}
to have $y_2\leq y\leq \min\big(1,\frac{\hat{\mu}_2}{\hat{\lambda}_2}\big)$, or $\max\big(1,\frac{\hat{\mu}_2}{\hat{\lambda}_2}\big)\leq y\leq y_3$. This means that for $1\leq x\leq\frac{\hat{\mu}_1}{\hat{\lambda}_1}$, $y_2 \leq Y_0(x)\leq 1$.
\qed
\begin{lem}\label{lem4-2}
 We present more properties about $Y_0(x)$ and $X_0(y)$ below:
\begin{description}
\item[(i)] If $\hat{\mu}_2>\hat{\lambda}_2$, then $0<Y_0(x)<1$ for $x\in\big(1,\frac{\hat{\mu}_1}{\hat{\lambda}_1}\big)$, and $1<Y_0(x)<\sqrt{\frac{\hat{\mu}_2}{\hat{\lambda}_2}}$ for $x\in\big(\frac{\hat{\mu}_1}{\hat{\lambda}_1},x_3\big)$. Specially, $Y_0\big(\frac{\hat{\mu}_1}{\hat{\lambda}_1}\big)=1$ and $Y_0(x_3)=\sqrt{\frac{\hat{\mu}_2}{\hat{\lambda}_2}}>1$.
\item[(ii)] If $\hat{\mu}_2<\hat{\lambda}_2$, then $0<Y_0(x)<1$ for $x\in(1,x_3)$. Also, $Y_0\big(\frac{\hat{\mu}_1}{\hat{\lambda}_1}\big)=\frac{\hat{\mu}_2}{\hat{\lambda}_2}<1$ and $Y_0(x_3)=\sqrt{\frac{\hat{\mu}_2}{\hat{\lambda}_2}}<1$.
\item[(iii)] If $\hat{\mu}_2=\hat{\lambda}_2$, then $x_3=\frac{\hat{\mu}_1}{\hat{\lambda}_1}$, $0<Y_0(x)<1$ for $x\in(1,\frac{\hat{\mu}_1}{\hat{\lambda}_1})$, and $Y_0(1)=Y_0\big(\frac{\hat{\mu}_1}{\hat{\lambda}_1}\big)=1$.
\end{description}
Similarly,
\begin{description}
\item[(i')] If $\hat{\mu}_2>\hat{\lambda}_2$, then $0<X_0(y)<1$ for $y\in\big(1,\frac{\hat{\mu}_2}{\hat{\lambda}_2}\big)$, and $1<X_0(y)<\sqrt{\frac{\hat{\mu}_1}{\hat{\lambda}_1}}$ for $y\in\big(\frac{\hat{\mu}_2}{\hat{\lambda}_2},y_3\big)$. Specially, $X_0\big(\frac{\hat{\mu}_2}{\hat{\lambda}_2}\big)=1$ and $X_0(y_3)=\sqrt{\frac{\hat{\mu}_1}{\hat{\lambda}_1}}>1$.
\item[(ii')] If $\hat{\mu}_2\leq\hat{\lambda}_2$, then $1<X_0(y)<\sqrt{\frac{\hat{\mu}_1}{\hat{\lambda}_1}}$ for $y\in(1,y_3)$. Also, $X_0(1)=1$ and $X_0(y_3)=\sqrt{\frac{\hat{\mu}_1}{\hat{\lambda}_1}}>1$.
\end{description}
\end{lem}

\proof
Based on equations (\ref{4-1})(\ref{4-3}), it is easy to know that $\hat{b}(x)<0$ for $x\in(1,x_3)$, so the branch $Y_0(x)$ should take $Y_-(x)=\frac{-\hat{b}(x)-\sqrt{\Delta(x)}}{2\hat{\lambda}_2 x}$. Solving the inequalities $Y_0(x)>1$ and $Y_0(x)<1$ with $x\in(1,x_3)$, we obtain, after some simple calculations, the results in (i)--(iii) of the lemma. (i') and (ii') can be proved in the same way.
\qed

\subsection{Poles of $\pi_1^{(1)}(x)$}

Since the censored random walk is a standard walk in the quarter plane, literature results can now be applied to the analysis of the dominant singularity of $\pi_1^{(1)}(x)$.
Therefore, besides the branch point $x_3$, given in (\ref{eqn:x3}), the other candidate for the dominant singularity is a pole of function $\pi_1^{(1)}(x)$.
In the following, we refine literature results, which lead to an explicit characterization of both the dominant pole and the regions for different exact tail asymptotic properties.

\begin{thm} {\bf (Theorem 4.4 in \cite{LZ2012})} \label{thm5-2}
 If $x_p$ is the pole of $\pi_1^{(1)}(x)$ with the smallest modulus in $(1, x_3]$, then $x_p$ is a zero of $h_1(x,Y_0(x))$ or $Y_0(x_p)$ is a zero of $h_2(X_0(y),y)$. In the latter case, $|Y_0(x_p)|>1$. On the other hand, if $x_p$ is the zero of $h_1(x,Y_0(x))$ or $Y_0(x_p)$ with $|Y_0(x_p)|>1$ is a zero of $h_2(X_0(y),y)$ with the smallest modulus in $(1, x_3]$, then $x_p$ is the pole of $\pi_1^{(1)}(x)$  in $(1, x_3]$. Moreover, $x_p$ is real. Parallel results can be easily stated for $\pi_2^{(1)}(y)$.
\end{thm}

For the retrial queue model with two input streams and two orbits studied in this paper, we show in the following that the pole of $\pi_1^{(1)}(x)$ (respectively $\pi_2^{(1)}(y)$) can only be the zero of $h_1(x,Y_0(x))$ (respectively $h_2(X_0(y),y)$).

First, we discuss properties of the pole of $\pi_1^{(1)}(x)$ in interval $(1,x_3]$. For convenience, let
$x^*$ be the unique zero in $(1,x_3]$ of the function $h_1(x,Y_0(x))$ if such a zero exists, otherwise let $x^*=+\infty$ (in this case, obviously $x^*$ can never be the dominant singularity since $x_3 < +\infty$).
Instead of directly considering the equation $h_1(x,Y_0(x))=0$, we consider the product of two functions $h_1(x,Y_0(x))$ and $h_1(x,Y_1(x))$, which is a polynomial:
\beqnn 
h_1(x,Y_0(x))h_1(x,Y_1(x))=\frac{\hat{\mu}_2}{\alpha (\lambda+\mu_1)^2}(x-1)g(x),
\eeqnn
where
\begin{equation*}
g(x)=\lambda\lambda_1(\lambda+\mu_1)x^2+\lambda\mu_1(\lambda+\mu_1-\mu)x-\mu\mu_1^2.
\end{equation*}
Since $\lambda(\lambda_1+\mu_1)<\mu\mu_1$, it is easily to check that $g(0)<0$ and $g(1)<0$. Hence $g(x)=0$ has one positive zero $x_+$ and one negative zero $x_-$. Especially, $x_+>1$. The expressions of the two zeros are given as
\begin{equation} \label{eqn:pole}
x_+=\frac{-\lambda\mu_1(\lambda+\mu_1-\mu)+\sqrt{[\lambda\mu_1(\lambda+\mu_1-\mu)]^2+4\lambda\lambda_1(\lambda+\mu_1)\mu\mu_1^2}}{2\lambda\lambda_1(\lambda+\mu_1)},
\end{equation}
\begin{equation*}
x_-=\frac{-\lambda\mu_1(\lambda+\mu_1-\mu)-\sqrt{[\lambda\mu_1(\lambda+\mu_1-\mu)]^2+4\lambda\lambda_1(\lambda+\mu_1)\mu\mu_1^2}}{2\lambda\lambda_1(\lambda+\mu_1)}.
\end{equation*}

Since either $p_{i,j}$ or $p_{i,j}^{(1)}$ is not X-shaped (refer to \cite{LZ2012} for details) in this censored random walk, Theorem~4.5 in Li and Zhao~\cite{LZ2012} guarantees that the candidate zero of $h_1(x,Y_0(x))$ can only be $x_+$. Solving $h_1(x_+,y)=0$, and then from \eqref{3-3} we get
\begin{equation} \label{5-14}
y=Y(x_+)=1-\frac{[\lambda_1(\lambda+\mu_1)x_+-\mu\mu_1](x_+-1)}{\lambda_2(\lambda+\mu_1)x_+},
\end{equation}
where $Y(x_+)$ is either $Y_0(x_+)$ or $Y_1(x_+)$.
On the other hand, $\frac{\mu\mu_1}{\lambda_1(\lambda+\mu_1)}>1$ and $g\big(\frac{\mu\mu_1}{\lambda_1(\lambda+\mu_1)}\big)=\frac{\lambda_2\mu\mu_1^{2}}{\lambda_1}>0$, hence, $1<x_+<\frac{\mu\mu_1}{\lambda_1(\lambda+\mu_1)}$. This means $Y(x_+)>1$.
Furthermore, to check whether or not $x_+$ is the pole of $\pi_1^{(1)}(x)$, we will carry out a discussion under the condition $\hat{\mu}_2>\hat{\lambda}_2$ and $\hat{\mu}_2\leq\hat{\lambda}_2$, respectively, in the following lemma.

\begin{lem}\label{lem5-3}
1. When $\hat{\mu}_2>\hat{\lambda}_2$, the value of $x^*$ depends on the value of $x_+$:
\begin{description}
\item[(a)] For $x_+\in\left(1,\frac{\hat{\mu}_1}{\hat{\lambda}_1}\right]$, we have $x^*=+\infty$;
\item[(b)] For $x_+\in\left(\frac{\hat{\mu}_1}{\hat{\lambda}_1},\min\big(x_3, \frac{\mu\mu_1}{\lambda_1(\lambda+\mu_1)}\big)\right)$, we have $x^*=x_+<x_3$ if $Y(x_+)<\sqrt{\frac{\hat{\mu}_2}{\hat{\lambda}_2}}$, and $x^*=+\infty$ otherwise;
\item[(c)] For $x_+=x_3<\frac{\mu\mu_1}{\lambda_1(\lambda+\mu_1)}$, we have $Y(x_+)=\sqrt{\frac{\hat{\mu}_2}{\hat{\lambda}_2}}$ and $x^*=x_+=x_3$;
\item[(d)] For $x_3<x_+<\frac{\mu\mu_1}{\lambda_1(\lambda+\mu_1)}$, we have $x^*=+\infty$.
\end{description}

2. When $\hat{\mu}_2\leq\hat{\lambda}_2$, we have $x^*=+\infty$.
\end{lem}

\proof
For the case $\hat{\mu}_2>\hat{\lambda}_2$, if $x_+\in\left(1,\frac{\hat{\mu}_1}{\hat{\lambda}_1}\right]$, it leads to $0<Y_0(x)<1$ from Lemma \ref{lem4-2}. For the case $\hat{\mu}_2\leq\hat{\lambda}_2$, if $x_+\in(1,x_3]$, it leads to $Y_0(x_+)\leq\sqrt{\frac{\hat{\mu}_2}{\hat{\lambda}_2}}\leq1$ from Lemma \ref{lem4-1}. The both cases contradict $Y(x_+)>1$. Hence, we can conclude that $h_1(x,Y_0(x))$ has no zero on $[1,+\infty)$ ($x_+$ is the zero of $h_1(x,Y_1(x))$). Therefore, $x^*=+\infty$.
Other conclusions are easy to make.
\qed

\begin{rem} It is worthwhile to notice that: (i) If there does not exist a pole in $(1,x_3]$, then $x_3$ is the dominant singularity of $\pi_1^{(1)}(x)$. Therefore, for the purpose of dominant singularity, we do not need to consider case 1(d) in Lemma~\ref{lem5-3}. (ii) The right-hand expression in (\ref{5-14}) can be either $Y_0(x_+)$ or $Y_1(x_+)$. (iii) It is possible that both $x_+$ and $x_-$ are zeros of $h_1(x,Y_1(x))$. In this case,  $h_1(x,Y_0(x))=0$ has no solution.
\end{rem}

Next, we show $h_2(X_0(y),y)$ has no zeros.
Based on Theorem \ref{thm5-2}, if the pole in $(1,x_3]$ of $\pi_1^{(1)}(x)$ is not $x^*$, then it is denoted by $\tilde{x}_1$. For convenience,
define $y^*$ to be the unique zero of $h_2(X_0(y),y)$ in $(1,y_3]$ if such a zero exists, otherwise let $y^*=+\infty$.
Following the same idea as above, we have
\beqnn 
h_2(X_0(y),y)h_2(X_1(y),y)=\frac{\hat{\mu}_1}{\alpha (\lambda+\mu_2)^2}(y-1)f(y),
\eeqnn
where
\begin{equation*} 
f(y)=\lambda\lambda_2(\lambda+\mu_2)y^2+\lambda\mu_2(\lambda+\mu_2-\mu)y-\mu\mu_2^2
\end{equation*}
has two zeros: $y_-<0$ and $y_+>1$. Solving $h_2(x,y_+)=0$, and then from (\ref{3-4}) we get
\beqnn 
x=X(y_+)=1-\frac{[\lambda_2(\lambda+\mu_2)y_+-\mu\mu_2](y_+-1)}{\lambda_1(\lambda+\mu_2)y_+},
\eeqnn
where $X(y_+)$ is either $X_0(y_+)$ or $X_1(y_+)$.

Using a similar argument, parallel results to Lemma \ref{lem5-3}-1 can be obtained. Since $\lambda_2(\lambda+\mu_2)<\mu\mu_2$ and $f\big(\frac{\mu\mu_2}{\lambda_2(\lambda+\mu_2)}\big)=\frac{\lambda_1\mu\mu_2^{2}}{\lambda_2}>0$, we have $1<y_+<\frac{\mu\mu_2}{\lambda_2(\lambda+\mu_2)}$. This leads to $X(y_+)>1$. Next, we claim that  $\tilde{x}_1$ cannot exist.

If $h_2(X_0(y),y)$ has a zero $y^*$ in $(1,y_3]$, then $y^*=y_+$ and $1<X(y_+)=X_0(y^*)\leq\sqrt{\frac{\hat{\mu}_1}{\hat{\lambda}_1}}$.
From Theorem 4.7 in \cite{LZ2012} we know $\tilde{x}_1=X_1(y^*)$. Hence, $\tilde{x}_1\in\left[\sqrt{\frac{\hat{\mu}_1}{\hat{\lambda}_1}},\frac{\hat{\mu}_1}{\hat{\lambda}_1}\right)$ and $0<Y_0(\tilde{x}_1)<1$ from Lemma \ref{lem4-2}. Obviously, it contradicts to that $Y_0(\tilde{x}_1)$ is a pole of $\pi_2^{(1)}(y)$. Therefore, $\tilde{x}_1$ cannot exist.

Based on the above discussion, we are ready to summarize the detailed properties on the location of the dominant singularity. For convenience, we introduce the following three conditions:
\begin{description}
\item[Condition~1.] $\hat{\mu}_2 > \hat{\lambda}_2$, $x_+\in\left(\frac{\hat{\mu}_1}{\hat{\lambda}_1},\min\big(x_3,\frac{\mu\mu_1}{\lambda_1(\lambda+\mu_1)}\big)\right)$ and $Y(x_+) < \sqrt{\frac{\hat{\mu}_2}{\hat{\lambda}_2}}$.
\item[Condition~2.] $\hat{\mu}_2 > \hat{\lambda}_2$ and $x_+ =x_3 \in \big(\frac{\hat{\mu}_1}{\hat{\lambda}_1},\frac{\mu\mu_1}{\lambda_1(\lambda+\mu_1)}\big)$.
\item[Condition~3.] One of the following three: (a) $\hat{\mu}_2 \leq \hat{\lambda}_2$; (b) $\hat{\mu}_2 > \hat{\lambda}_2$ and $x_+ \in (1, \frac{\hat{\mu}_1}{\hat{\lambda}_1}]$; and
(c) $\hat{\mu}_2 > \hat{\lambda}_2$, $x_+\in\left(\frac{\hat{\mu}_1}{\hat{\lambda}_1},\min\big(x_3,\frac{\mu\mu_1}{\lambda_1(\lambda+\mu_1)}\big)\right)$, and $Y(x_+) \geq \sqrt{\frac{\hat{\mu}_2}{\hat{\lambda}_2}}$.
\end{description}

\begin{lem}\label{lem6-1}
\begin{description}
\item[Case 1:] Under Condition~1, the dominant singularity $x_{dom}=x^*=x_+ < x_3$, which is a pole.
\item[Case 2:] Under Condition~2, the dominant singularity $x_{dom}=x_3=x^*=x_+$, which is both a branch point and a pole.
\item[Case 3:] Under Condition~3, the dominant singularity $x_{dom}=x_3 < x^*=+\infty$, which is a branch point.
\end{description}
\end{lem}

\begin{rem}
One should notice that the above lemma is a refinement of the literature result for a general random walk in the quarter plane. It provides explicit conditions (in terms of system parameters), under which the dominant singularity $x_{dom}$ of $\pi_1^{(1)}(x)$ (also explicitly expressed) is either $x_{dom} = x_3$ or $x_{dom}=x_+$, since all the branch point $x_3$, the pole $x_+$ and $Y(x_+)$ are explicitly expressed in (\ref{eqn:x3}), (\ref{eqn:pole}) and (\ref{5-14}), respectively.
\end{rem}

\section{Asymptotic properties of $\pi_1^{(1)}(x)$ at its dominant singularity}

Once again, in this section, we only provide detailed analysis for the function $\pi_1^{(1)}(x)$. Due to symmetry, parallel results for $\pi_2^{(1)}(y)$ can be easily stated and similarly proved. In the previous section, we proved that either $x_3$ or $x_+$ is the dominant singularity of $\pi_1^{(1)}(x)$. In this section,
   we prove (in Theorem~\ref{thm6-2}) that there exist three types of detailed asymptotic properties as $x$ approaches to the dominant singularity $x_{dom}$ of $\pi_1^{(1)}(x)$, depending on $x_{dom}=x_+<x_3$ or $x_{dom}=x_3<x_+$ or $x_{dom}=x_+=x_3$ respectively.

For simplicity in the following discussion, especially for the case of $x_{dom}=x_3$, we write
\beqlb\label{6-1}
Y_0(x)=p(x)+q(x)\sqrt{1-\frac{x}{x_{dom}}}\;,
\eeqlb
\beqnn 
h_1(x,Y_0(x))=p_1(x)+q_1(x)\sqrt{1-\frac{x}{x_{dom}}}\;,
\eeqnn
\beqnn 
Y_0(x_{dom})-Y_0(x)=\bigg(1-\frac{x}{x_{dom}}\bigg)p^*(x)-q(x)\sqrt{1-\frac{x}{x_{dom}}}\;,
\eeqnn
\beqnn 
h_1(x,Y_0(x))-h_1(x_{dom},Y_0(x_{dom}))=\bigg(1-\frac{x}{x_{dom}}\bigg)p_1^*(x)+q_1(x)\sqrt{1-\frac{x}{x_{dom}}}\;,
\eeqnn
where
\beqnn
p(x)=-\frac{\hat{b}(x)}{2\hat{\lambda}_2x}, ~~ q(x)=-\frac{1}{2\hat{\lambda}_2x}\sqrt{\frac{\Delta(x)}{1-x/x_{dom}}},~~\text{if}\;~~ x_{dom}=x_3,
\eeqnn
\beqnn
p_1(x)=-\frac{\hat{b}(x)a_1(x)}{2\hat{\lambda}_2x}+b_1(x), ~~ q_1(x)=a_1(x)q(x),
\eeqnn
\beqnn
p^*(x)=\frac{(p(x_{dom})-p(x))x_{dom}}{x_{dom}-x}~~ \text{and}\;~~p_1^*(x)=\frac{(p_1(x)-p_1(x_{dom}))x_{dom}}{x_{dom}-x}.
\eeqnn

\begin{thm}\label{thm6-2}
The behaviour of $\pi_1^{(1)}(x)$ at the dominant singularity is given as
\begin{description}
\item[(i)] If $x_{dom}=x^*=x_+<x_3$, then
\beqnn
\displaystyle\lim_{x\rightarrow x_+}\bigg(1-\frac{x}{x_+}\bigg)\pi_1^{(1)}(x)=C_{1,0},
\eeqnn
where
\beqnn
C_{1,0}=\frac{(\lambda+\mu_1)\sqrt{\Delta(x_+)}\left[h_2(x_+,Y_0(x_+))\pi_2^{(1)}(Y_0(x_+))+h_0(x_+,Y_0(x_+))\pi_{0,0}(1)\right]}{\hat{\mu}_2\lambda\lambda_1x_+(x_+-1)(x_+-x_-)}.
\eeqnn
\item[(ii)] If $x_{dom}=x_3=x^*=x_+$, then
\beqnn
\displaystyle\lim_{x\rightarrow x_{dom}}\sqrt{1-x/x_{dom}}\pi_1^{(1)}(x)=C_{2,0},
\eeqnn
where
\beqnn
C_{2,0}=\frac{2 }{\lambda_1}\times\frac{h_2(x_{dom},Y_0(x_{dom}))\pi_2^{(1)}(Y_0(x_{dom}))+h_0(x_{dom},Y_0(x_{dom}))\pi_{0,0}(1)}{\sqrt{x_{dom}(x_{dom}-x_1)(x_{dom}-x_2)(x_4-x_{dom})}}.
\eeqnn
\item[(iii)] If $x_{dom}=x_3<x^*=+\infty$, then
\beqnn
\displaystyle\lim_{x\rightarrow x_3}\sqrt{1-x/x_3}\pi_1^{'(1)}(x)=C_{3,0},
\eeqnn
where $\pi_1^{'(1)}(x)$ is the derivative of $\pi_1^{(1)}(x)$ and
\beqnn
C_{3,0}=-\frac{q(x_3)}{2x_3}\frac{d}{dy}\left[\frac{h_2(x_3,y)\pi_2^{(1)}(y)+h_0(x_3,y)\pi_{0,0}(1)}{h_1(x_3,y)}\right]{\bigg|_{y=Y_0(x_3)}}.
\eeqnn
\end{description}
\end{thm}

\proof (i) If $x_{dom}=x^*=x_+<x_3$, then $x_{dom}$ is a simple pole of $\pi_1^{(1)}(x)$. Based on the analysis in \cite{LZ2012}, we can rewrite
\beqnn
\pi_1^{(1)}(x)&=&-\frac{\left[h_2(x,Y_0(x))\pi_2^{(1)}(Y_0(x))+h_0(x,Y_0(x))\pi_{0,0}(1)\right]h_1(x,Y_1(x))}{h_1(x,Y_0(x))h_1(x,Y_1(x))}
\\&=&-\frac{\left[h_2(x,Y_0(x))\pi_2^{(1)}(Y_0(x))+h_0(x,Y_0(x))\pi_{0,0}(1)\right]h_1(x,Y_1(x))}{\frac{\hat{\mu}_2}{\alpha (\lambda+\mu_1)^2}(x-1)g(x)}
\\&=&-\frac{\left[h_2(x,Y_0(x))\pi_2^{(1)}(Y_0(x))+h_0(x,Y_0(x))\pi_{0,0}(1)\right]h_1(x,Y_1(x))}{\frac{\hat{\mu}_2}{\alpha (\lambda+\mu_1)}(x-1)\lambda\lambda_1(x-x_-)(x-x_+)}.
\eeqnn
It follows that
\beqnn
\displaystyle\lim_{x\rightarrow x_+}\bigg(1-\frac{x}{x_+}\bigg)\pi_1^{(1)}(x)=C_{1,0}.
\eeqnn

(ii) If $x_{dom}=x_3=x^*=x_+$, then $h_1(x_{dom},Y_0(x_{dom}))=0$. In this case,  we can rewrite $\pi_1^{(1)}(x)$ as
\beqnn
\pi_1^{(1)}(x)=\frac{-h_2(x,Y_0(x))\pi_2^{(1)}(Y_0(x))-h_0(x,Y_0(x))\pi_{0,0}(1)}{\sqrt{1-x/x_{dom}}\left[\sqrt{1-x/x_{dom}}p_1^*(x)+q_1(x)\right]}.
\eeqnn
It follows that
\beqnn
\displaystyle\lim_{x\rightarrow x_{dom}}\sqrt{1-x/x_{dom}}\pi_1^{(1)}(x)=\frac{h_2(x_{dom},Y_0(x_{dom}))\pi_2^{(1)}(Y_0(x_{dom}))+h_0(x_{dom},Y_0(x_{dom}))\pi_{0,0}(1)}{-a_1(x_{dom})q(x_{dom})}=C_{2,0}.
\eeqnn

(iii) If $x_{dom}=x_3<x^*$, let
\beqnn
T(x,y)=\frac{-h_2(x,Y_0(x))\pi_2^{(1)}(Y_0(x))-h_0(x,Y_0(x))\pi_{0,0}(1)}{h_1(x,Y_0(x))}.
\eeqnn
Then the derivative of $\pi_1^{(1)}(x)$ is given by
\beqnn
    \pi_1^{'(1)}(x)=\frac{\partial T}{\partial x}+\frac{\partial T}{\partial y}\frac{dY_0(x)}{dx}
\eeqnn
with
\beqnn
\frac{dY_0(x)}{dx}=p'(x)+q'(x)\sqrt{1-x/x_{dom}}-\frac{q(x)}{2x_{dom}\sqrt{1-x/x_{dom}}},
\eeqnn
where $p(x)$ and $q(x)$ are defined in equation \eqref{6-1}. Since it is obvious that $\displaystyle\lim_{x\rightarrow x_3}\sqrt{1-x/x_3}$$\frac{dY_0(x)}{dx}=-\frac{q(x_3)}{2x_3}$, $\displaystyle\lim_{x\rightarrow x_3}\sqrt{1-x/x_3}$$\frac{\partial T}{\partial x}=0$ and $\frac{\partial T}{\partial y}$ is continuous at $(x_3,Y_0(x_3))$, so,
\beqnn
\displaystyle\lim_{x\rightarrow x_{dom}}\sqrt{1-x/x_3}\pi_1^{'(1)}(x)&=&-\frac{q(x_3)}{2x_3}\frac{\partial T}{\partial y}\bigg|_{(x_3,Y_0(x_3))}=C_{3,0}.
\eeqnn
\qed

\section{Tail asymptotic properties in stationary probabilities}

Exact tail asymptotic properties in stationary probabilities are obtained directly from the corresponding asymptotic properties of the unknown generating function by applying the following Tauberian-like theorem. This theorem is originated from
Bender~\cite{B1974}, and  more complete versions can be found in Flajolet and Sedgewick~\cite{FS2009}, which include the following theorem as a special case.

\begin{thm}\label{thm6-1}
{\bf (Tauberian-like theorem for single singularity)}
Let $A(z)=\sum_{n\geq 0}a_n z^n$ be analytic at zero with $R$ the radius of convergence.
Suppose that $R$ is a singularity of $A(z)$ that can be continued to a $\Delta$-domain at $R$. If for a real number
$\beta \notin \{0,-1,-2, \ldots \}$,
\beqnn
\lim_{z\to R}(1-z/R)^{\beta}A(z)=g,
\eeqnn
where $g$ is a  non-zero constant. Then,
\beqnn
a_n \sim \frac{g}{\Gamma(\beta)}n^{\beta-1}R^{-n},
\eeqnn
where $\Gamma(\beta)$ is the value of Gamma function at $\beta$, and $a_n \sim b_n$ is equivalent to $\lim_n a_n/b_n =1$.
\end{thm}

The Tauberian-like theorem claims that the tail behaviour in the sequence of the coefficients in the Taylor expansion of the analytic function corresponds to the asymptotic property of the function at its dominant singularity. In the following subsections, we show how to apply Theorem~\ref{thm6-1} to characterize the tail behaviour in the joint probabilities $\pi_{m,n}(k)$
for a fixed number $n$ of customers in orbit 2. Specifically, in subsection~1, we provide
a characterization for tail asympotics, when the server is busy, in the sequence of: (1) boundary probabilities $\pi_{m,0}(1)$; (2) marginal probabilities $\pi_m^{(1)}=\sum_{n=1}^{\infty}\pi_{m,n}(1)$; (3) joint probabilities $\pi_{m,n}(1)$ for a fixed $n > 0$ (along the direction of queue one). While in subsection~2, when the server is idle, we provide a characterization for tail asympotics in $\pi_{m,n}(0)$ for a fixed $n$ and for the marginal distribution $\pi_m^{(0)}=\sum_{n=1}^{\infty}\pi_{m,n}(0)$.

\begin{rem}
By symmetry, tail behaviour in $\pi_{m,n}(1)$ and $\pi_{m,n}(0)$ for a fixed number $m$ of customers in orbit 1 (and also in the marginal distributions for the second queue length when the server is busy and idle, respectively) can be easily stated and similarly proved.
\end{rem}

\subsection{Exact tail asymptotics when the server is busy}

First, we consider the sequence $\pi_{m,0}(1)$ of the boundary probabilities. When the second queue is empty and the server is busy, the exact tail asymptotic behaviour of the stationary probability sequence $\pi_{m,0}(1)$ along the increasing direction of the first queue is a direct consequence of the characterization of the asymptotic property for the function $\pi_1^{(1)}(x)$ in Theorem~\ref{thm6-2} and the Tauberian-like theorem (Theorem~\ref{thm6-1}).

\begin{thm} \label{theorem:6.2}
 For a stable retrial queue with two input streams and two orbits studied in this paper, when $m$ is large, we have three types of tail asymptotic properties for the boundary probabilities $\pi_{m,0}(1)$:
\begin{description}
\item[Type~1: (Exact geometric decay)] Under Condition~1,
\beqnn 
\pi_{m,0}(1)\sim C_{1,0}\left(\frac{1}{x_+}\right)^{m-1}, \quad m\geq1;
\eeqnn
\item[Type~2: (Geometric decay with prefactor $m^{-1/2}$)] Under Condition~2,
\beqnn 
\pi_{m,0}(1)\sim \frac{C_{2,0}}{\sqrt{\pi}}m^{-\frac{1}{2}}\left(\frac{1}{x_{dom}}\right)^{m-1}, \quad m\geq1;
\eeqnn
\item[Type~3: (Geometric decay with prefactor $m^{-3/2}$)] Under Condition~3,
\beqnn 
\pi_{m,0}(1)\sim \frac{C_{3,0}}{\sqrt{\pi}}m^{-\frac{3}{2}}\left(\frac{1}{x_3}\right)^{m-2}, \quad m\geq1.
\eeqnn
\end{description}
Here, constants $C_{i,0}$ ($i=1,2,3$) are given in Theorem \ref{thm6-2}.
\end{thm}

\begin{rem} One may notice that in Type~3, the power of the decay rate is $m-2$ instead of $m-1$ since the Tauberian-like theorem is applied to the derivative of the function.
\end{rem}

For characterizing the asymptotic behaviour of the marginal probability $\pi_m^{(1)}=\sum_{n=1}^{\infty}\pi_{m,n}(1)$, we compute $\pi^{(1)}(x,1)$,
\beqnn
\pi^{(1)}(x,1)&=&\frac{h_1(x,1)\pi_1^{(1)}(x)+h_2(x,1)\pi_2^{(1)}(1)+h_0(x,1)\pi_{0,0}(1)}{-h(x,1)}
\\&=&-\frac{\frac{1}{(\lambda+\mu_1)}\left[\lambda_1(\lambda+\mu_1)x-\hat{\mu}_1\right]\pi_1^{(1)}(x)+\lambda_1\pi_2^{(1)}(1)+\lambda_1\pi_{0,0}(1)}{\lambda_1x-\frac{\hat{\mu}_1}{\alpha}}
\\&=&\frac{\alpha}{\lambda+\mu_1}\frac{\left[\lambda_1(\lambda+\mu_1)x-\hat{\mu}_1\right]\pi_1^{(1)}(x)+\lambda_1(\lambda+\mu_1)(\pi_2^{(1)}(1)+\pi_{0,0}(1))}{\hat{\mu}_1(1-\frac{\hat{\lambda}_1}{\hat{\mu}_1}x)}.
\eeqnn

If $\hat{\lambda}_2\neq\hat{\mu}_2$, it follows from (\ref{4-2}) that we have $\hat{\mu}_1/\hat{\lambda}_1<x_3$. Therefore, from Lemma~\ref{lem6-1} we can claim that $1<\hat{\mu}_1/\hat{\lambda}_1<\min(x^*,x_3)$ is always true. Obviously, $\hat{\mu}_1/\hat{\lambda}_1$ is the dominant singularity of $\pi^{(1)}(x,1)$, which is a simple pole.
If $\hat{\lambda}_2=\hat{\mu}_2$, then from Lemma~\ref{lem4-2}-(iii), we have $x_3=\hat{\mu}_1/\hat{\lambda}_1$. Again, according to Lemma~\ref{lem6-1}, the dominant singularity of $\pi_1^{(1)}(x)$ is $x_3=\hat{\mu}_1/\hat{\lambda}_1<x^*=+\infty$. Notice that $\lim_{x\rightarrow x_3} \pi_1^{(1)}(x)$ is finite. Therefore, the Tauberian-like theorem can be still applied.
\begin{thm}\label{theorem:6.3}
(i)
\[
    \displaystyle\lim_{x\rightarrow\hat{\mu}_1/\hat{\lambda}_1}\bigg(1-\frac{x}{\hat{\mu}_1/\hat{\lambda}_1}\bigg)\pi^{(1)}(x,1)=C_m,
\]
where
\beqlb\label{eq6-8}
C_m=
-\frac{\mu_2}{\lambda+\mu_1}\pi_1^{(1)}(\hat{\mu}_1/\hat{\lambda}_1)+\frac{\hat{\lambda}_1}{\hat{\mu}_1}\big(\pi_2^{(1)}(1)+\pi_{0,0}(1)\big);
\eeqlb
and (ii)
The marginal probabilities $\pi_m^{(1)}$ has an exact geometric decay with decay rate equal to $x_{dom}=\hat{\mu}_1/\hat{\lambda}_1$:
\beqnn
    \pi_m^{(1)}\sim C_m\left(\frac{\hat{\lambda}_1}{\hat{\mu}_1}\right)^{m-1}.
\eeqnn
\end{thm}

\begin{rem}
It should be noticed that one may consider $\sum_{n=0}^{\infty} (\pi_{m,n}(1)+\pi_{m,n}(0))$ the usual marginal distribution of the first queue. Its tail asymptotic property can be easily obtained since the property for $\pi_m^{(1)}$ and $\pi_{m,0}(1)$ have been studied, and the property for $\pi_m^{(0)}=\sum_{n=1}^{\infty}\pi_{m,n}(0)$ and $\pi_{m,0}(0)$ can be similarly obtained.
\end{rem}

Next, the exact tail asymptotic behaviour for joint probabilities can be obtained from the recursive relationship of the generating functions $\varphi_n(x)$,
defined by
\[
\varphi_n(x)=\sum_{m=1}^{\infty}\pi_{m,n}(1)x^{m-1}, \quad n\geq0.
\]
It is clear that $\varphi_0(x)=\pi_1^{(1)}(x)$. From the balance equations of the censored random walk, we can obtain
\begin{eqnarray}
\label{6-8}
c(x)\varphi_1(x)+b_1(x)\varphi_0(x)&=&a_0^*(x),  \\
\label{6-9}
c(x)\varphi_2(x)+b(x)\varphi_1(x)+a_1(x)\varphi_0(x)&=&a_1^*(x), \\
\label{6-10}
c(x)\varphi_{n+1}(x)+b(x)\varphi_n(x)+a(x)\varphi_{n-1}(x)&=&a_n^*(x), \quad n\geq2,
\end{eqnarray}
where
\begin{eqnarray*}
a_0^*(x)&=&-c_2(x)\pi_{0,1}-b_0(x)\pi_{0,0}, \\
a_1^*(x)&=&-c_2(x)\pi_{0,2}-b_2(x)\pi_{0,1}-a_0(x)\pi_{0,0}, \\
a_n^*(x)&=&-c_2(x)\pi_{0,n+1}-b_2(x)\pi_{0,n}-a_2(x)\pi_{0,n-1}, \quad n\geq2.
\end{eqnarray*}
Rewrite (\ref{6-10}) as
\beqnn
\varphi_{n+1}(x)=\frac{-b(x)\varphi_n(x)-a(x)\varphi_{n-1}(x)+a_n^*(x)}{c(x)}, \quad n\geq2,
\eeqnn
and note that $c(x)=p_{0,-1}x$. Hence, we established the fact that $\varphi_n(x)$ has the same singularities as $\varphi_0(x)$ since that the zero of $c(x)$ is not a pole of $\varphi_n(x)$ for all $n\geq0$.

By adopting Theorem 7.1 and Lemma 7.2 in \cite{LZ2012} directly, we define
\beqnn 
A_i(x_{dom})=-\frac{b_1(x_{dom})}{c(x_{dom})}C_{i,0}, ~~ i=1,2,3 \quad \text{and} \quad B_3(x_3)=-\frac{p_1(x_3)}{c(x_3)}C_{3,0},
\eeqnn
then we can conclude the results in the following theorem.
\begin{thm} \label{theorem:6.4}
Corresponding to the three types in Theorem \ref{thm6-2}, when $m$ is large, we have the following tail asymptotic properties for the joint probabilities $\pi_{m,n}(1)$ for a fixed $n$:
\begin{description}
\item[Type~1: (Exact geometric decay)]
\beqnn 
\pi_{m,n}(1)\sim A_1(x_+)\bigg(\frac{1}{Y_1(x_+)}\bigg)^{n-1}\left(\frac{1}{x_+}\right)^{m-1}, \quad n\geq1;
\eeqnn
\item[Type~2: (Geometric decay with prefactor $m^{-1/2}$)]
\beqnn 
\pi_{m,n}(1)\sim \frac{A_2(x_{dom})}{\sqrt{\pi}}\bigg(\frac{1}{Y_1(x_{dom})}\bigg)^{n-1}m^{-\frac{1}{2}}\left(\frac{1}{x_{dom}}\right)^{m-1}, \quad n\geq1;
\eeqnn
\item[Type~3: (Geometric decay with prefactor $m^{-3/2}$)]
\beqnn 
\pi_{m,n}(1)\sim \frac{[A_3(x_3)+(n-1)B_3(x_3)]}{\sqrt{\pi}}\bigg(\frac{1}{Y_1(x_3)}\bigg)^{n-1}m^{-\frac{3}{2}}\left(\frac{1}{x_3}\right)^{m-2}, \quad n\geq1.
\eeqnn
\end{description}
\end{thm}

\subsection{Exact tail asymptotics when the server is idle}

Having known the exact tail asymptotic properties of the boundary, marginal and joint distributions for $I(t)=1$ (or the server is busy), we can now study the tail asymptotic properties for $I(t)=0$ (or the server is idle) based on the relationship given in (\ref{2-4}).

Setting $y=0$ in \eqref{2-4} leads to
\beqlb\label{eq6-5}
(\lambda+\mu_1)P^{(0)}(x,0)=\mu P^{(1)}(x,0)+\mu_1P^{(0)}(0,0),
\eeqlb
which means that $P^{(0)}(x,0)$ and $P^{(1)}(x,0)$ have the same asymptotic property.

Similarly, setting $y=1$ in \eqref{2-4} leads to
\beqlb\label{eq6-6}
\alpha P^{(0)}(x,1)=\mu P^{(1)}(x,1)+\mu_2 P^{(0)}(x,0)+\mu_1 P^{(0)}(0,1).
\eeqlb
Substituting \eqref{eq6-5} into \eqref{eq6-6} gives
\beqnn
\alpha P^{(0)}(x,1)=\mu P^{(1)}(x,1)+\frac{\mu\mu_2}{\lambda+\mu_1} P^{(1)}(x,0)+\frac{\mu_1\mu_2}{\lambda+\mu_1}P^{(0)}(0,0) + \mu_1 P^{(0)}(0,1).
\eeqnn
Since the asymptotic property at the dominant singularity of $P^{(0)}(x,1)$ is dominated by the asymptotic property of the function $\mu P^{(1)}(x,1)$, $P^{(0)}(x,1)$ and $P^{(1)}(x,1)$
have the same asymptotic property. Based on the above, we have the following conclusion:

\begin{thm} \label{theorem:7.1}
Assume that the retrial queue with two input streams and two orbits is stable.
\begin{description}
\item[(i)] For large $m$, we have three types of tail asymptotic properties for the boundary probabilities $\pi_{m,0}(0)$ correspondingly.
\begin{description}
\item[Type~1: (Exact geometric decay)]
\beqnn 
\pi_{m,0}(0)\sim \frac{\mu}{\lambda+\mu_1}C_{1,0}\left(\frac{1}{x_+}\right)^{m-1}, \quad m\geq1;
\eeqnn
\item[Type~2: (Geometric decay with prefactor $m^{-1/2}$)]
\beqnn 
\pi_{m,0}(0)\sim \frac{\mu}{\lambda+\mu_1}\frac{C_{2,0}}{\sqrt{\pi}}m^{-\frac{1}{2}}\left(\frac{1}{x_{dom}}\right)^{m-1} \quad m\geq1;
\eeqnn
\item[Type~3: (Geometric decay with prefactor $m^{-3/2}$)]
\beqnn 
\pi_{m,0}(0)\sim \frac{\mu}{\lambda+\mu_1}\frac{C_{3,0}}{\sqrt{\pi}}m^{-\frac{3}{2}}\left(\frac{1}{x_3}\right)^{m-2} \quad m\geq1.
\eeqnn
\end{description}
Here, constants $C_{i,0}$ $(i=1,2,3)$ are given in Theorem \ref{thm6-2}.

\item[(ii)] The tail asymptotic property of the marginal distribution $\pi_m^{(0)}=\sum_{n=1}^{\infty}\pi_{m,n}(0)$ is determined by
\beqnn
\pi_m^{(0)}\sim \frac{\mu}{\alpha}C_m\left(\frac{\hat{\lambda}_1}{\hat{\mu}_1}\right)^{m-1},
\eeqnn
where $C_m$ is provided by \eqref{eq6-8}.
\end{description}
\end{thm}

We finally characterize the tail asymptotic behaviour for the joint probabilities $\pi_{m,n}(0)$ for a fixed $n > 0$.
Define the generating function
\beqnn
G_n^{(k)}(x)=\sum_{m=0}^{\infty}\pi_{m,n}(k)x^{m}, \quad k=0,1 ~~ n\geq1.
\eeqnn
Referring to equation (14) in \cite{ANY2014}, we have
\beqnn
\alpha G_n^{(0)}(x)-\mu G_n^{(1)}(x)=\mu_1\pi_{0,n}(0),
\eeqnn
which obviously leads to the following theorem.

\begin{thm} \label{theorem:7.2}
Corresponding to the three types in Theorem \ref{thm6-2}, when $m$ is large, we have the following tail asymptotic properties of the joint probabilities $\pi_{m,n}(0)$ for a fixed $n$:
\begin{description}
\item[Type~1: (Exact geometric decay)] Under Condition~1,
\beqnn 
\pi_{m,n}(0)\sim \frac{\mu}{\alpha}A_1(x_+)\bigg(\frac{1}{Y_1(x_+)}\bigg)^{n-1}\left(\frac{1}{x_+}\right)^{m-1}, \quad n\geq1;
\eeqnn
\item[Type~2: (Geometric decay with prefactor $m^{-1/2}$)] Under Condition~2,
\beqnn 
\pi_{m,n}(0)\sim \frac{\mu}{\alpha}\frac{A_2(x_{dom})}{\sqrt{\pi}}\bigg(\frac{1}{Y_1(x_{dom})}\bigg)^{n-1}m^{-\frac{1}{2}}\left(\frac{1}{x_{dom}}\right)^{m-1}, \quad n\geq1;
\eeqnn
\item[Type~3: (Geometric decay with prefactor $m^{-3/2}$)] Under Condition~3,
\beqnn 
\pi_{m,n}(0)\sim \frac{\mu}{\alpha}\frac{[A_3(x_3)+(n-1)B_3(x_3)]}{\sqrt{\pi}}\bigg(\frac{1}{Y_1(x_3)}\bigg)^{n-1}m^{-\frac{3}{2}}\left(\frac{1}{x_3}\right)^{m-2}, \quad n\geq1.
\eeqnn
\end{description}
\end{thm}

\section{Concluding remarks}

In this paper, we considered the exact tail asymptotic behaviours of a retrial queue with two input streams and two orbits. Partitioned according to the two states of the server, this model is formulated as a random walk in the quarter plane whose transition probabilities are modulated by a two-state Markov chain (idle or busy). Our work is a revisit of the same model studied in \cite{ANY2014}. While in \cite{ANY2014}, the study is based on the solution to a BVP, we employed a different method, the kernel method. The advantage of using this method mainly relies on the fact that there is no need to have a full determination of the unknown generating function. Instead, we only need the location of the dominant singularity of the unknown function and the asymptotic property at its dominant singularity. By this method, tail asymptotic properties in stationary probabilities for the model are obtained when the first queue size is large.
Due to symmetry, it is not difficult to state and (similarly) prove parallel exact tail asymptotic properties when the second queue size is large. In addition, exact tail asymptotic results for other probability sequences formed from the joint stationary probabilities can also be considered. For example, we can consider the total number of customers in the system as follows: let
\[
    \pi_T=\sum_{\substack{m,n: \\ m+n=T}}\pi_{m,n}
\]
 and we compute $\pi^{(1)}(x,x)$,
 according to \eqref{3-15}:
\beqnn
\pi^{(1)}(x,x)=-\frac{\big(\lambda x-\frac{\hat{\mu}_1}{\lambda+\mu_1}\big)\pi^{(1)}_1(x)+\big(\lambda x-\frac{\hat{\mu}_2}{\lambda+\mu_2}\big)\pi^{(1)}_2(x)+\lambda\pi_{0,0}(1)}{x(\lambda x-\frac{\hat{\mu}_1+\hat{\mu}_2}{\alpha})}.
\eeqnn
Then, the dominant singularity is determined by comparing $x=(\hat{\mu}_1+\hat{\mu}_2)/\hat{\lambda}$ to the dominant singularities of $\pi^{(1)}_1(x)$ and $\pi^{(1)}_2(x)$, and therefore the asymptotic property at its dominant singularity is determined. The exact tail asymptotic property is a consequence of the Tauberian-like theorem.

This paper used a censored chain to convert the matrix-form fundamental form into a usual (scalar) fundamental form. It is not always feasible to do this conversion since explicit expressions might not exist for the censored chain. A general method is to solve the matrix-form fundamental form to have a relationship between generating functions for different states of the modulated chain. For example, the censored chain to the idle state does not have an explicit expression for its transition matrix. However, in terms of the relationship in \eqref{2-4} and \eqref{2-5} obtained by solving the matrix-form fundamental form, we can have the following functional equation:
\[
R(x,y)P^{(0)}(x,y)=A(x,y)P^{(0)}(x,0)+B(x,y)P^{(0)}(0,y),\quad |x|\leq1,|y|\leq1,
\]
with
\begin{align*}
R(x,y)&=\hat{\lambda}_1(1-x)xy+\hat{\lambda}_2(1-y)xy-\hat{\mu}_1(1-x)y-\hat{\mu}_2(1-y)x,
\\A(x,y)&=\left[(1-y)(\lambda_2 y-\mu)+\lambda_1(1-x)y\right]\mu_2 x,
\\B(x,y)&=\left[(1-x)(\lambda_1 x-\mu)+\lambda_2(1-y)x\right]\mu_1 y,
\end{align*}
which is equivalent to:
\begin{equation*}
R(x,y)\pi^{(0)}(x,y)=\frac{A(x,y)-R(x,y)}{y}\pi_1^{(0)}(x)+\frac{B(x,y)-R(x,y)}{x}\pi_2^{(0)}(y)+\frac{A(x,y)+B(x,y)-R(x,y)}{xy}\pi_{0,0}(0).
\end{equation*}
After some calculations, the above equation also can be written as
\[
-\hat{h}^{(0)}(x,y)\pi^{(0)}(x,y)=\hat{h}_1^{(0)}(x,y)\pi_1^{(0)}(x)+\hat{h}_2^{(0)}(x,y)\pi_2^{(0)}(y)+\hat{h}_0^{(0)}(x,y)\pi_{0,0}(0),
\]
where
\begin{align*}
\hat{h}^{(0)}(x,y)&=[\hat{\lambda}_1x+\hat{\lambda}_2y+\hat{\mu}_1x^{-1}+\hat{\mu}_2y^{-1}-(\hat{\lambda}+\hat{\mu}_1+\hat{\mu}_2)]xy,
\\\hat{h}_1^{(0)}(x,y)&=[\lambda_1(\lambda+\mu_1)x+\lambda_2(\lambda+\mu_1)y+\hat{\mu}_1x^{-1}-\lambda(\lambda+\mu_1)-\hat{\mu}_1]x,
\\\hat{h}_2^{(0)}(x,y)&=[\lambda_1(\lambda+\mu_2)x+\lambda_2(\lambda+\mu_2)y+\hat{\mu}_2y^{-1}-\lambda(\lambda+\mu_2)-\hat{\mu}_2]y,
\\\hat{h}_0^{(0)}(x,y)&=\lambda\lambda_1x+\lambda\lambda_2y-\lambda^2.
\end{align*}
The above functional equation is the fundamental form corresponding a random walk defined by
\beqnn
\hat{p}_{1,0}=\hat{\lambda}_1, \quad \hat{p}_{0,1}=\hat{\lambda}_2, \quad \hat{p}_{-1,0}=\hat{\mu}_1, \quad \hat{p}_{0,-1}=\hat{\mu}_2,  \quad \hat{p}_{0,0}=1-(\hat{\lambda}+\hat{\mu}_1+\hat{\mu}_2),
\eeqnn
\beqnn
\hat{p}_{1,0}^{(1)}=\lambda_1(\lambda+\mu_1), \quad \hat{p}_{0,1}^{(1)}=\lambda_2(\lambda+\mu_1), \quad \hat{p}_{-1,0}^{(1)}=\hat{\mu}_1,  \quad \hat{p}_{0,0}^{(1)}=1-[\lambda(\lambda+\mu_1)+\hat{\mu}_1],
\eeqnn
\beqnn
\hat{p}_{1,0}^{(2)}=\lambda_1(\lambda+\mu_2), \quad \hat{p}_{0,1}^{(2)}=\lambda_2(\lambda+\mu_2), \quad \hat{p}_{0,-1}^{(2)}=\hat{\mu}_2, \quad \hat{p}_{0,0}^{(2)}=1-[\lambda(\lambda+\mu_2)+\hat{\mu}_2],
\eeqnn
\beqnn
\hat{p}_{1,0}^{(0)}=\lambda\lambda_1, \quad \hat{p}_{0,1}^{(0)}=\lambda\lambda_2,  \quad \hat{p}_{0,0}^{(0)}=1-\lambda^2.
\eeqnn
We now can apply the kernel method to the resulting fundamental form to obtain exact tail asymptotic properties for probabilities with an idle server.

Finally, we emphasize that this work serves as an illustration of how the kernel method can be applied to random walks modulated by a finite-state Markov chain, which has a similar structure property to that the retrial queue model possesses.

\vspace*{5mm}
\noindent{\bf Acknowledgments:} This work was done during the visit of the first two authors
to the School of Mathematics and Statistics,  Carleton University (Ottawa, Canada), who acknowledge the support provided by the School. The first author also thanks the China Scholarship Council for supporting her visit to Carleton University through a scholarship. In addition, this work was supported in partial by the National Natural Science Foundation of China (11271373), and by the Natural Sciences and Engineering Research Council of Canada (NSERC). All authors thank the comments/suggestions made by two anonymous reviewers, which significantly improved the quality of the paper.

\end{document}